\input epsf.sty
\input graphicx
\input  amssym.tex
\def\sqr#1#2{{\vcenter{\hrule height.#2pt              
     \hbox{\vrule width.#2pt height#1pt\kern#1pt
     \vrule width.#2pt}
     \hrule height.#2pt}}}
\def\square{\mathchoice\sqr{5.5}4\sqr{5.0}4\sqr{4.8}3\sqr{4.8}3}
\def\qed{\hskip4pt plus1fill\ $\square$\par\medbreak}

\def\P{{\bf P}}


\centerline{\bf Linear Fractional Recurrences:}

\centerline{\bf  Periodicities and Integrability}
\bigskip
\centerline{Eric Bedford and Kyounghee Kim}

\bigskip\bigskip

\noindent{\bf \S0.  Introduction.}  Let $k\ge3$, and let $\alpha=(\alpha_0,\dots,\alpha_k)$, $\beta=(\beta_0,\dots,\beta_k)$ be $(k+1)$-tuples of complex numbers.  We consider a $k$-step linear fractional recurrence
$$x_{n+k+1} = {\alpha_0 +  \alpha_1 x_{n+1}+ \cdots + \alpha_k x_{n+k} \over \beta_0 + \beta_1 x_{n+1} +\cdots + \beta_k x_{n+k}} \eqno(0.1)$$
Given a $k$-tuple $(x_1,\dots,x_k)$, the relation (0.1) generates a sequence $\{x_j, j\ge1\}$ as long as the denominator does not vanish.  The question has been raised (see [GL] and [CL]) to find the $\alpha$ and $\beta$ for which (0.1) is periodic.  By ``periodic'' we mean that the sequence $\{x_j,j\ge0\}$ is periodic for every starting point $(x_1,\dots,x_k)$.  There are a number of works in the literature that have considered this question under the hypothesis that all numbers are positive.  Here we consider it natural to examine this question over the field of complex numbers.  

The case $k=2$ was considered in [BK1] for general $\alpha$ and  $\beta$, and it was shown that the only possible nontrivial periods are 6, 5, 8, 12, 18, and 30.  (Here, ``nontrivial'' means that the map cannot be reduced to a simpler map, e.g.\ linear or 1-dimensional.)  McMullen [M] observed that these periods are the orders of certain Coxeter groups and showed that these maps represent the corresponding Coxeter elements.  The case of dimension 3 is determined in [BK2]:  the only  possible nontrivial periods are 8 and 12.  The 3-step, period 8 maps had been found previously; there are two essentially different maps, one is in [L], and the other is in [CsLa].  Here we show that the period 12 corresponds to a phenomenon that holds for $k$-step recurrences for all $k$:
\proclaim Theorem 0.1.  For each $k$, there are $k$ different recurrences of the form (0.1) with $\alpha$, $\beta$ as in (5.3), which have period $4k$.

Our approach is similar to that of [BK1,2]: we consider (0.1) in terms of the associated birational map 
$$f_{\alpha,\beta}(x_1,\dots,x_k) = \left ( x_2,\dots,x_k,{\alpha_0 +  \alpha_1 x_1+ \cdots + \alpha_k x_k \over \beta_0 + \beta_1 x_1 +\cdots + \beta_k x_k} \right) \eqno(0.2)  $$
of $k$-dimensional space.  We may consider $f_{\alpha,\beta}$ as a birational map of complex projective space ${\bf P}^k$, as well as any space $X$ which is birationally equivalent to ${\bf P}^k$.  For a rational map $f$ of $X$, there are well-defined pull-back maps $f^*$ on the cohomology groups $H^{p,q}(X)$, as well as on the Picard group $Pic(X)$.  
We may define a notion of growth by 
$$\delta(f):=\lim_{n\to\infty}||(f^n)^*||^{1\over n}.$$
Here we work on $Pic(X)$ (or $H^{1,1}$), where $\delta(f)$ is equivalent to degree growth.    To determine $\delta(f)$, we replace ${\bf P}^k$ with a space $X$ with the property that passage from $f$ to $f^*$ is compatible with iteration.  Specifically, we ``regularize'' the map $f$ in the sense that we replace ${\bf P}^k$ by an  $X$ such that $(f^n)^*=(f^*)^n$ holds on $Pic(X)$; so in this case we obtain $\delta(f)$ as the modulus of the largest eigenvalue of $f^*$. 
  The way we find our space $X$ is to analyze the ``singular'' behavior of $f$, by which we mean the behavior that prevents $f_{\alpha,\beta}$ from being a diffeomorphism.  Namely, there are hypersurfaces $E$ with the property that either $f(E)$ or $f^{-1}(E)$ has codimension $>1$.  Such a hypersurface is called {\it exceptional}.    The existence/nonexistence of exceptional hypersurfaces depends on the choice of representative $X$ for $f$, and the regularity of $f_X$ is determined by the behavior of the orbits of exceptional hypersurfaces.  

In \S2 we show that for generic $\alpha$ and $\beta$ we have $\delta(f_{\alpha,\beta})= \Delta_k>1$.   In particular, we conclude that a generic $f_{\alpha,\beta}$ is not periodic.  In order to prove Theorem 0.1, we find a space $X$ for which $f_X^*$ is periodic, and we use this to conclude that $f$ is periodic.   Although our map $f_X$ is not an automorphism,  de Fernex and Ein [dFE] have shown that since $f$ is periodic there will exist a space $Z$ so that the induced map $f_Z$ is biholomorphic.  

Next we consider the mappings 
$$h(x_1,\dots,x_k) = \left(x_2,\dots,x_k,{a + x_2 + x_3+\cdots + x_k\over x_1}\right), \eqno(0.3)$$
which have been discussed in several places, often under the name of ``Lyness map'' because of its origin in [L].  Except in two exceptional cases, these maps are not periodic, but they exhibit an integrability which has been studied by several authors: ([KLR], [KL], [Z], [CGM1--3], [GBM], [HKY], [GKI]).   Applying our analysis to $h$ we construct a rather different regularization and obtain:
\proclaim Theorem 0.2.  If $k>3$, or if  $a\ne1$, the degree of $h^n$ is quadratic in $n$.

In dimension $k=2$, there is a strong connection (see [DF] and [G]) between polynomial degree growth and integrability.  Namely, if  $g$ is a birational surface map, then  linear degree growth corresponds to preserving  a  rational fibration; and quadratic degree growth corresponds to preserving an elliptic fibration.  In \S6 we discuss the structure of rational functions that are invariant under $f$ which were found in [GKI].

\bigskip
\noindent{\bf \S1.  Birational maps.} Let us recall a few notions from algebraic geometry that we will use.  The reader is referred to [B] for further details.   A rational map of projective space ${\bf P}^k$  is given by a $k+1$-tuple of homogeneous polynomials $f=[f_0:f_1:\cdots:f_k]$ of a common degree $d={\rm deg}(f_0)=\cdots={\rm deg}(f_k)$.  We refer to $d$ as the degree of $f$.  Without loss of generality we assume that the $f_j$'s have no common factor, so the degree is well defined.  The indeterminacy locus is ${\cal I}(f) = \{x:f_0(x)=\cdots=f_k(x)=0\}$.   Since the $f_j$'s have no common factor, ${\cal I}(f)$ always has codimension at least 2.  The map $f$ is holomorphic at all points of  ${\bf P}^k-{\cal I}(f)$ but cannot be extended to be continuous at any point of ${\cal I}$.  If $V$ is a variety for which no irreducible component is contained in ${\cal I}$ then the {\it strict transform} $f(V)$ is defined as the closure of $f(V-{\cal I})$.  The strict transform is again a variety.  We say that a map $f$ is {\it dominant} if its image contains an open set.  Given two rational maps $f$ and $g$, there is a rational map $f\circ g$, and $f\circ g$ is equal to $f(g(x))$ for all $x\notin {\cal I}(g)$ such that $g(x)\notin{\cal I}(f)$.   The map $f$ is said to be {\it birational} if there is a rational map  $g$ such that $f\circ g$ and $g\circ f$ are both the identity.  

If $f$ is a rational map, we say that a subvariety $E$ is {\it exceptional} if $E\notin{\cal I}$, and the dimension of $f(E)$ is strictly less than the dimension of $E$.  We let ${\cal E}$ denote the set of exceptional hypersurfaces of  $f$.  We will say that $f$ is a {\it pseudo-automorphism} if there is no exceptional hypersurface. 

We will define manifolds by the procedure of blowing up.  If  $p\in X$ is a point, then the {\it blowup} of $X$ at $p$ is given by a new manifold $Y$ with a holomorphic projection $\pi:Y \to X$ such that $\pi:Y-\pi^{-1}p\to X-p$ is biholomorphic, and $\pi^{-1}p$ is equivalent to ${\bf P}^{k-1}$.  Similarly, if  $S$ is a smooth submanifold of $X$, we may define a blowup of  $X$ along the center $S$.  Given a blowup $\pi:Y\to X$, the preimage $\pi^{-1}S$ of $S$ under $\pi$ will be called the exceptional blowup fiber.  If $f:X\dasharrow X$ is a rational map, then there is an induced rational map $f_Y:=\pi^{-1}\circ f_X\circ \pi$ on $Y$.

We refer to a variety of pure codimension 1 as a hypersurface.  We say that two hypersurfaces $S_1$ and $S_2$ are {\it linearly equivalent} if there is a rational function $r$ such that the divisor of $\{r=0\}$ is  $S_1-S_2$.  By $Pic(X)$ we denote the set of all divisors modulo linear equivalence.  The spaces $X$ that we will deal with all arise from ${\bf P}^k$ by blowups, so that in fact $Pic(X)$ is isomorphic to the Dolbeault cohomology group $H^{1,1}(X)$.  If $\pi:Y\to X$ is a blowup along a smooth center $S$, and if ${\cal S}$ denotes the exceptional blowup fiber over $S$, then $Pic(Y)$ is generated by $Pic(X)$, together with the class of ${\cal S}$.

Suppose that $f:X\dasharrow X$ is a rational map.  Given a hypersurface $S=\{p=0\}$ we define the pullback  $f^*S$ as the closure of $\{p(f)=0\}-{\cal I}(f)$.    This gives a well defined linear map $f^*$ of $Pic(X)$.   We say that $f$ is {\it 1-regular} if $(f^n)^*=(f^*)^n$ holds on $Pic(X)$.  The (first) dynamical degree of $f_X$ is defined by the growth of the iterates on $Pic(X)$:
$$\delta(f_X):=\lim_{n\to\infty}||(f^n_X)^*|_{Pic(X)}||^{1\over n}$$
This is independent of the choice of norm $||\cdot||$.  And since $\pi:X\to{\bf P}^k$ is holomorphic, we have $(f_X)^n=(\pi^{-1}\circ f\circ \pi)^n=\pi^{-1}\circ f^n\circ\pi = (f^n)_X$, so $\delta(f)=\delta(f_X)$.

\bigskip
\noindent{\bf \S 2. Linear Fractional Recurrences }  We may interpret equation $(0.2)$ as a rational map $f : \P^k \dasharrow \P^k$ by writing it in homogeneous coordinates as 
%

$$ f[x_0 : \cdots : x_k] = [ x_0 \beta \cdot x : x_2 \beta \cdot x : \dots : x_k \beta \cdot x : x_0 \alpha\cdot x] \eqno{(2.1)}$$
where $\alpha \cdot x = \alpha_0 x_0 + \alpha_1 x_1 + \cdots \alpha_k x_k$. Let us set $ B= (- \alpha_1, 0, \dots, 0, \beta_1)$, $\alpha'=( \alpha_0, \alpha_2, \dots, \alpha_k, 0)$ and $\beta' =(\beta_0, \beta_2, \dots, \beta_k, 0)$. The inverse of $f$ is given by the map 
$$f^{-1}[x_0 : \cdots : x_k] = [ x_0 B\cdot x : x_0 \alpha' \cdot x - x_k \beta' \cdot x : x_1 B\cdot x : \cdots : x_{k-1} B \cdot x] \eqno{(2.2)}$$ 
We assume that 
$$\eqalign{ & (\beta_1,\beta_2,\dots, \beta_k) \ne (0,0,\dots,0),\cr & \alpha {\rm \ is\ not\ a\ multiple\ of \ }\beta, {\rm\ \  and}\phantom{qwerqwuytre} \cr
& (\alpha_i,\beta_i) \ne (0,0) {\rm \ for\ } i=1\   {\rm \ and \ some\ } 1 < i \le k.\cr}\eqno{(2.3)}$$
If $(\alpha_1,\beta_1)=(0,0)$ then $f$ does not depend on $x_1$ and thus $f$ can be realized as a $k-1$ step recurrence relation. If $(\alpha_i,\beta_i)=(0,0)$ for all $i=2, \dots,k$ then $f^k$ is an essentially $1$-dimensional mapping.

Let us set $\gamma= \beta_1 \alpha- \alpha_1 \beta$ and $C= \beta_1 \alpha' - \alpha_1 \beta'$. For $0 \le i \le k$ we use notation $\Sigma_i = \{ x_i = 0 \}$ and  $e_i = [ 0 : \cdots: 0:1:0:\cdots:0]$, the point whose $i$-th coordinate is nonzero and everything else is zero. We also use $\Sigma_\beta  = \{ \beta\cdot x = 0\}, \Sigma_\gamma = \{ \gamma \cdot x = 0\}, \Sigma_B = \{ B \cdot x = 0\},$ and $\Sigma_C= \{ C \cdot x = 0 \}$.  To indicate the intersection we combine their subscripts, for example 
$ \Sigma_{0\beta}= \{ x_0=\beta\cdot x=0\}$ and $ \Sigma_{01}= \{ x_0= x_1=0\}$. 
The Jacobian of $f$ is a constant multiple of $x_0 (\beta \cdot x)^{k-1} (\gamma \cdot x).$ The Jacobian vanishes on three hypersurfaces $\Sigma_0, \Sigma_\beta, \Sigma_\gamma$; these hypersurfaces are exceptional and are mapped to the lower dimensional linear subspaces:
$$f(\Sigma_0 ) = \Sigma_{0B}=\Sigma_{0k},\quad f(\Sigma_\beta) = e_k,\quad f(\Sigma_\gamma) = \Sigma_{BC}.$$
The Jacobian of $f^{-1}$ is a constant multiple of $x_0 (B\cdot x)^{k-1}(C \cdot x)$ and we have 
$$f^{-1} : \Sigma_0 \mapsto \Sigma_{0\beta},\ \ \ \Sigma_B \mapsto e_1,\ \ \ {\rm and} \ \ \ \Sigma_C \mapsto \Sigma_{\beta\gamma}.$$ The indeterminacy locus of $f$ is ${\cal I}^+ = \{ e_1, \Sigma_{0\beta}, \Sigma_{\beta\gamma}\}$ and $f^{-1}$ is ${\cal I}^+ = \{ e_k, \Sigma_{0B}, \Sigma_{BC}\}$. 

\medskip
Let us consider the maps which satisfy (2.3) and the following:
$$ \beta_1 \ne 0 \qquad {\rm and } \qquad \beta_1\alpha_j-\alpha_1\beta_j \ne 0\ \ {\rm for\  all\ }\  j=2, \dots, k\eqno{(2.4)}$$
For every choice of parameters $\alpha,\beta$ satisfying (2.3--4),  we have 
$$f: \Sigma_0 \mapsto \Sigma_{0k}\mapsto \Sigma_{0\,k-1\,k}\mapsto \cdots \mapsto \Sigma_{0\,3\,\dots\,k} \mapsto e_1 \rightsquigarrow \Sigma_B \eqno{(2.5)}$$
We first modify the orbit of $\Sigma_0$. Let  $\pi:Y \to {\bf P}^k$ be the complex manifold obtained by blowing up $e_1$ and then $\Sigma_{0\,3\,\dots\,k} $ and continuing successively until we reach $\Sigma_{0\,k}$.  That is, we let $\Sigma_{3\, \dots\, k}$ denote its strict transform in the space obtained by blowing up $e_1$, and then we blow up along the center $\Sigma_{3\, \dots\, k}$, etc.  We let $E_1=\pi^{-1} e_1$ denote the exceptional fiber over $e_1$, and let  ${\cal S}_{0,j}$ denote the exceptional fiber over $\Sigma_{0\,j\,\dots\,k}$ for all $j \ge 3$.  Let us set $\alpha^{(j)}=( 0,\alpha_1, \dots, \alpha_j,0,\dots,0)$, and $\beta^{(j)}=(0,\beta_1, \dots, \beta_j,0,\dots,0)$. For $3 \le j \le k$ we use local coordinates near ${\cal S}_{0\,j}$
$$ \pi_{sj}: (s, x_2, \dots, x_{j-1}, \xi_j, \dots, \xi_{k} ) _{s_j} \mapsto [ s:1: x_2: \cdots:x_{j-1}: s \xi_j:\cdots:s \xi_k] \in \P^k$$
and for the neighborhood of the exceptional divisor $E_1$ we use 
$$ \pi_{e_1} : ( s, \xi_2, \dots, \xi_k)_{e_1} \mapsto [s:1:s \xi_2: \cdots: \xi_k] \in \P^k.$$
Working with the induced birational map $f_Y:=\pi^{-1}\circ f\circ\pi$ we have 
$$f_Y:\Sigma_0\ni [0:x_1:x_2:\cdots, x_k] \mapsto \left(0,x_2, \dots, x_{k}, {\alpha^{(k)}\cdot x \over \beta^{(k)} \cdot x}\right)_{s_k} \in {\cal S}_{0\, k}$$
Similarly we have for all $j=2, \dots, k-1$
$$f_Y:{\cal S}_{0,j+1}\ni (0,x_2, \dots, x_{j}, \xi_{j+1}\dots ,\xi_k)_{s_{j+1}}
\mapsto\left(0, x_3/x_2, \dots, x_j/x_2, \xi_{j+1}, \dots, \xi_k,{\alpha^{(k)}\cdot x \over \beta^{(k)} \cdot x}\right)_{s_j} \in{\cal S}_{0\, j}    $$
For the points of $E_1$, we have
$$f_Y:E_1\ni (0,\xi_2, \dots, \xi_k) \mapsto [ \beta_1: \beta_1 \xi_2: \cdots: \beta_1 \xi_k: \alpha_1] \in \Sigma_B.$$
By condition $(2.4)$, we see that $\alpha^{(j)}$ is not a constant multiple of $\beta^{(j)}$ for all $2 \le j\le k$. It follows that 
\proclaim Lemma 2.1. The map $f_Y$ is a local diffeomorphism at generic points of $\Sigma_0, E_1$, and $S_{0,j+1}$ for $2 \le j \le k-1$.  

Since the induced map $f_Y$ is a local diffeomorphism at points of $\Sigma_0 \cup E_1 \cup \bigcup_{j = 3}^k {\cal S}_{0,j}$, $\Sigma_0$ and all the exceptional (blowup) divisors $E_1$ and ${\cal S}_{0,j}$ for $j=3,\dots,k$ are not exceptional for $f_Y$. Thus the exceptional set for $f_Y$ consists of two divisors:  ${\cal E}^+_Y = \{ \Sigma_\beta, \Sigma_\gamma\}$. The indeterminacy locus for $f_Y$ is ${\cal I}^+_Y = \{ \Sigma_{0\beta},\Sigma_{\beta \gamma}\}$. For the inverse map $f_Y^{-1}$ we have ${\cal E}^-_Y = \{ \Sigma_0, \Sigma_C\}$ and ${\cal I}^-_Y = \{ \{e_k\},\Sigma_{BC}\}$.

 Let us consider the ordered basis of $Pic(Y)$ : $H$, $E_1$, ${\cal S}_{0,3}, \dots,$ ${\cal S}_{0,k}$. Using the discussion above, we have 
\proclaim Lemma 2.2. With the given ordered basis the action of $f_Y^*$ on $Pic(Y)$ is given by
$$\eqalign{f_Y^*\ :\  & E_1 \mapsto {\cal S}_{0,3}\mapsto \cdots \mapsto {\cal S}_{0,k}\cr
& {\cal S}_{0,k}\mapsto \{\Sigma_0\} = H-E_1-{\cal S}_{0,3} - \cdots-{\cal S}_{0,k}\cr
& H \mapsto 2H - E_1.}$$
Under the same ordered basis the action of $f_Y^{-1}$  is:
$$\eqalign{{f_Y^{-1}}^*\ :\  &{\cal S}_{0,k} \mapsto {\cal S}_{0,k-1}\mapsto \cdots \mapsto {\cal S}_{0,3}\mapsto E_1\cr
& E_1\mapsto \{\Sigma_B\} = H-E_1-{\cal S}_{0,3} - \cdots-{\cal S}_{0,k}\cr
& H \mapsto 2H - E_1-{\cal S}_{0,3} - \cdots-{\cal S}_{0,k}.}$$

Now let us consider the following condition on exceptional hypersurfaces
$$\eqalign{ 
&f_Y^n {\cal E} \not\subset \Sigma_{\beta\gamma}\cup \Sigma_{0\beta} \quad {\rm \ for\ all\ \ } n \ge0 \qquad \qquad {\cal E} =\Sigma_\beta,\Sigma_\gamma  \cr 
& f_Y^{-n} {\cal E} \not\subset \Sigma_{BC} \cup \{e_k\}\quad {\rm \ for\ all\ \ } n \ge0 \qquad \qquad  {\cal E} =\Sigma_0,\Sigma_C \cr} \eqno{(2.6)}$$
When $(2.6)$ holds, $f_Y$ is regular in the sense that $(f_Y^*)^n=(f_Y^n)^*$ holds on $Pic(Y)$, and similarly for $f_Y^{-1}$.  Thus the dynamical degree $\delta(f)$ is the spectral radius of $f^*_Y$ acting on $Pic(Y)$. 

\proclaim Theorem 2.3.  For generic parameters,  the dynamical degrees satisfy the properties:  
\item{(i)} $\delta(f)=\Delta_k^+$ is the largest root of $x^k-(x^k-1)/(x-1)$, 
\item{(ii)} $\delta(f^{-1}) = \Delta_k^-$ is the largest root of $x^k-x^{k-1}+(-1)^{k-1}.$ 
\vskip0pt\noindent Furthermore we have
$$\lim_{k\to \infty} \Delta^+_k = 2, \quad {\rm and}\quad \lim_{k\to \infty} \Delta^-_k = 1.$$

\noindent{\it Proof.}   It is clear that for a generic map $f$, the critical triangle is nondegenerate.  Next we claim that (2.6) holds for a generic map.  Since (2.6) defines the complement of countably many varieties inside ${\bf C}^{2k+2}$, it suffices to show that there is one parameter value $(\alpha,\beta)$ for which (2.6) holds.  For each $k\ge 3$ let us consider $\alpha=(0,1,0,\dots,0)$ and $\beta=(0,1,\dots,2-k)$. 
In this case $\gamma = (0,0,-1,\dots,k-2)$ and $B=(1,0, \dots, 0,-1)$ and $C=(0,-1,\dots,-1,k-2,0)$. It follows that 
$$f_Y \Sigma_\gamma = \Sigma_{B\,C} \ni [ 1: 1: \cdots:1] \quad {\rm and} \quad f_Y[1: \dots:1] = [1:\dots:1]$$ 
We also have $f^{k+1}_Y:\Sigma_\beta =[1:0:\cdots:0:1]$ and
$$f_Y : [1:0:\cdots:0:1] \mapsto [1:0:\cdots:0:1:0] \mapsto\cdots\mapsto [1:1:0:\cdots:0] \mapsto[1:0:\cdots:0:1]$$
Thus $\Sigma_\beta$ is pre-periodic and there exists a fixed point in $\Sigma_\gamma \setminus \Sigma_\beta$. It follows that this mapping satisfies the the first half of (2.6).  The second part of condition (2.6) follows similarly.

The matrix representation of $f_Y^*$ and ${f_Y^{-1}}^*$ is given by $k\times k$ matrices:
$$f_Y^* = \pmatrix{ 2&0& \cdots&0  &1\cr-1&0& & &-1\cr 0 & 1& & &  \vdots\cr &0&\ddots& &\vdots\cr 0& & 0& 1&-1},\ \ {f_Y^{-1}}^* = \pmatrix{ 2&1& &  &0\cr-1&-1&1 & &0\cr \vdots& \vdots& &\ddots & \cr \vdots& \vdots& & &1\cr -1&-1& 0& 0&0}$$
The dynamical degrees $\Delta_k^+$ and $\Delta_k^-$ are given by the spectral radius of the above matrix representations.  Now (2.6) implies that $f_Y$ is regular, so the spectral radius of $f_Y^*$ gives the dynamical degree.  The polynomials in statements $(i)$ and $(ii)$ are the characteristic polynomials of these matrices.  \qed

\medskip
\proclaim Corollary 2.4. For every $f$ of the form (0.2), we have $\delta(f)\le \Delta_k^+$ and $\delta(f^{-1})\le \Delta_k^-$.  

\noindent{\it Proof.}   Thus we have $\delta(f_{\alpha,\beta})=\Delta_k^+$ and $\delta(f_{\alpha,\beta}^{-1})=\Delta_k^-$ for generic $\alpha$ and $\beta$.  The general inequality follows because the function $(\alpha,\beta)\mapsto \delta(f_{\alpha,\beta})$ is lower semicontinuous.  \qed

\bigskip

\noindent{\bf \S 3. Non-Periodicity }  We call the set of exceptional hypersurfaces $\{ \Sigma_0,\Sigma_\beta, \Sigma_\gamma\}$ the critical triangle. When these three hypersurfaces are distinct, we say the critical triangle is nondegenerate.

\proclaim Lemma 3.1. If $\beta_1 \ne 0$ and $\beta_1\alpha_j-\alpha_1\beta_j = 0$ for all $j\ne0$ then there is unique exceptional hypersurface. If $\beta_1=0$, then there are only two distinct exceptional hypersurfaces. Otherwise the critical triangle is nondegenerate. 

\noindent{\it Proof.} Using the condition $(1.3)$ we see that $\Sigma_0 \ne \Sigma_\beta$ and $\Sigma_B \ne \Sigma_C$. It follows that $f$ and $f^{-1}$ has at least two distinct exceptional hypersurfaces. If $\beta_1=0$, we have $\gamma= -\alpha_1 \beta$ and $B=(-\alpha_1,0,\dots,0)$. It follows that $\Sigma_\beta = \Sigma_\gamma$ and $\Sigma_0 = \Sigma_B$. If $\beta_1 \ne 0$ and $\beta_1\alpha_j-\alpha_1\beta_j = 0$ for all $j\ne0$, then $\gamma= C=\beta_1 (\alpha_0,0,\dots,0)$ and therefore $\Sigma_0= \Sigma_\gamma=\Sigma_C$.\qed

For $j=2, \dots, k$ let us consider the codimension $k-j+2$ linear subspaces 
$$ {\cal L}_j := \Sigma_{0\beta} \cap \bigcap_{\ell = j}^k \{\beta_1 x_{k-\ell+1} + \beta_2 x_{k-\ell+2} +\cdots \beta_j x_{k-\ell+k} = 0 \}.$$

\proclaim Lemma 3.2. Let $j_*$ be the largest integer such that $\beta_j \ne 0$. If $j_* >1$ then $e_k \not \in {\cal L}_{j_*}$ and $\Sigma_0$ is pre-fixed under $f^{-1}$: 
$$f^{-(k-j_*+1)}\Sigma_0 = {\cal L}_{j_*}, \quad f^{-1} {\cal L}_{j_*} = {\cal L}_{j_*}.$$
In case $j_* =1$ we have 
$$f^{-1}: \Sigma_0 \mapsto \Sigma_{0\,1}\mapsto \Sigma_{0\,1\,2}\mapsto \cdots\mapsto e_k \rightsquigarrow \Sigma_\beta.$$

\noindent{\it Proof.} In case $j_*>1$, since $\beta_{j_*} \ne 0$ it follows that $e_k \not \in {\cal L}_{j_*}$. The second part is the immediate consequence of $(1.5)$.To show that ${\cal L}_{j_*}$ is fixed under $f^{-1}$, first notice that ${\cal L}_{j_*}$ has codimension $2+k-j_*$. A generic point  $p\in {\cal L}_{j_*}$ can be written in terms of $x_{k-j_*+2}, \dots, x_k$ and $x_{k-j_*+1} = -(\beta_2 x_{k-j_*+2}+ \cdots + \beta_{j_*} x_k) $.  It follows that the image of this point $f^{-1}p$ is $[0:y_1: \cdots: y_k]$ where $y_i = x_{i-1}$ for $i \ge 2$ and there for codimension of $f^{-1} {\cal L}_{j_*}$ is  $2+k-j_*$.  When $j_* =1$, ${\cal L}_j = \Sigma_{0\,1\,\cdots\,k+1-j}$ for $j = 1, \dots, k$.\qed

If the mapping $f$ is periodic with period $p$ then $f^{-1}$ is also periodic and for every hypersurface $H$ in $\P^k$ $f^pH = H$ and therefore the codimension of $f^pH$ has to be equal to 1. Thus we have

\proclaim Corollary 3.3. If $\beta_j \ne 0$ for some $j\ge 2$, then $f$ is not periodic.

\bigskip

\noindent{\bf \S 4. Critical Case }   We say $f$ is {\it critical} if $\beta_j = 0$ for all $j>1$ and the critical triangle is non-degenerate. Using $(1.4)$ we may also set $\alpha_k=1$, and by Lemma 2.1, we may assume that
$$\alpha= (\alpha_0,0,\alpha_2, \dots, \alpha_{k-1},1)\ \ {\rm and\ \ }\beta=(\beta_0,1,0,\dots,0), \alpha_2 \cdots \alpha_{k-1} \ne 0. \eqno{(4.1)}$$
Let us consider the involution $\tau[x_0:x_1:\cdots :x_k]=[x_0:x_k:\cdots:x_1]$ gotten by interchanging the variables  $x_j \leftrightarrow x_{k-j+1}$, $1\le j\le k$.  We see that $f$ is reversible in the sense that  $f^{-1}=\tau\circ f\circ \tau$.
If $(4.1)$ holds, we have $$\gamma=\beta_1\alpha-\alpha_1\beta= \alpha, \quad B= (0,\dots,0,1),\quad {\rm and\ \ \ }C=\beta_1\alpha'-\alpha_1\beta' = \alpha'.$$
When the mapping is critical, we use the conjugacy by $\tau$ and apply Lemma 2.2 to $f^{-1}$ to obtain:
$$f_Y :\Sigma_\beta \mapsto e_k \rightsquigarrow \Sigma_{0\,1\,\dots\,k-2}\rightsquigarrow \Sigma_{0\,1\,\dots
,k-3}\rightsquigarrow\cdots \rightsquigarrow \Sigma_{0\,1} = \Sigma_{0\beta} \rightsquigarrow \Sigma_0\eqno{(4.2)}$$
where $e_k$, $\Sigma_{0\,1\,\dots\, k-2}$, $\dots$ , $\Sigma_{0\,1}$ are the strict transforms in $Y$ of the corresponding linear subspaces in ${\bf P}^k$. 

Let us consider a complex manifold $\pi_X:X \to Y$ obtained by a successive blowing up the sets $e_k$, $\Sigma_{0\,1\,\dots\, k-2}$, $\dots$, $\Sigma_{0\,1}$. We denote the exceptional divisors over  $e_k$, $\Sigma_{0\,1\,\dots\, k-2}$, $\dots$, $\Sigma_{0\,1}$ be $E_k$, ${\cal P}_{0,k-2}$, $\dots$, ${\cal P}_{0,1}$.  We will see that the induced maps on the blowup fibers are dominant:
$$f_X:\Sigma_\beta \to E_k\mapsto {\cal P}_{0,k-2}\mapsto\cdots\mapsto {\cal P}_{0,1} \mapsto \Sigma_0\mapsto {\cal S}_{0,k} \mapsto \cdots\mapsto {\cal S}_{0,3}\mapsto E_1\mapsto \Sigma_B.\eqno{(4.3)}$$
To see the mapping is dominant, we also work with local coordinates. For near $E_k$
$$\pi_{e_k} : (s, \xi_1, \dots, \xi_{k-1})_{e_k} \mapsto[s: s \xi_1: \cdots: s \xi_{k-1}:1]$$
and for $1 \le j \le k-2$, we use the following local coordinates for the neighborhood of ${\cal P}_{0,j}$ such that $\{s=0\} = {\cal P}_{0,j}$:
$$\pi_{pj} : ( s, \xi_1, \dots, \xi_j, x_{j+1}, \dots, x_{k-1})_{p_j} \mapsto [ s: s \xi_1: \cdots: s \xi_j : x_{j+1} : \cdots: x_{k-1}:1]$$
Then the induced birational map $f_X$ acts on $\Sigma_\beta$ and each exceptional divisors as follows:
$$\eqalign{f_X: &\Sigma_\beta \ni [x_0: -\beta_0 x_0: x_2: \cdots: x_k] \mapsto (0, x_2/x_0, \dots, x_k/x_0)_{e_k} \in E_k\cr
& E_k \ni (0, \xi_1, \dots, \xi_k)_{e_k} \mapsto (0, \xi_2, \dots, \xi_{k-1}, \beta_0 + \xi_1)_{p_{k-2}} \in {\cal P}_{0,k-2}}$$
For all $2 \le j\le k-2$
$$f_X: {\cal P}_{0,j} \ni(0,\xi_1, \dots, \xi_j, x_{j+1}, \dots, x_{k-1})_{p_j} \mapsto (0, \xi_2, \dots, \xi_j,{\beta\cdot \xi\, x_{j-1} \over \alpha \cdot x}, \dots, {\beta\cdot \xi\, x_{k-1} \over \alpha \cdot x})_{p_{j-1}} \in {\cal P}_{0,j-1}$$
where $\beta\cdot \xi = \beta_0 + \xi_1$ and $\alpha \cdot x = \alpha_{j+1} x_{j+1} + \cdots + \alpha_{k-1} x_{k-1} + \alpha_k$.
And we also have 
$$ f_X: {\cal P}_{0,1} \ni(0,\xi_1, x_2, \dots, x_{k-1})_{p_1} \mapsto [0: x_2 \beta\cdot \xi: \cdots: x_{k-1}\beta \cdot \xi: \beta \cdot \xi: \alpha_2 x_2 + \cdots + \alpha_{k-1} x_{k-1} + \alpha_k] \in \Sigma_0$$
On the other hand the induced map $f_X^{-1}$ acts as followings:
$$\eqalign{ f_X^{-1} : &\Sigma_0 \ni [0: x_1: \cdots: x_k] \mapsto (0, (\alpha'\cdot x- \beta_0 x_k) /x_k, x_1/x_{k-1}, \dots, x_{k-2}/ x_{k-1})_{p_1} \in {\cal P}_{0,1}\cr &E_k \ni (0,\xi_1, \dots, \xi_{k-1} )_{e_k} \mapsto [ 1: -\beta_0 : \xi_1: \cdots: \xi_{k-1}] \in \Sigma_\beta}$$
Hence $f_X$ is a local diffeomorphisms at points on $\Sigma_\beta \cup E_k \cup \bigcup^{k-2}_{j=1}{\cal P}_{0,j}$ and $f_X^{-1}$ is a local diffeomorphisms at points on $\Sigma_0 \cup E_k \cup \bigcup^{k-2}_{j=1}{\cal P}_{0,j}$. It follows that
\proclaim Lemma 4.1. In the critical case, the induced map $f_X$ has only one exceptional hypersurface $\Sigma_\gamma$; $\Sigma_{\beta\gamma}$ is the only component of the indeterminacy locus ${\cal I}(f_X)$ which blows  up to a hypersurface. 

\medskip
\proclaim Lemma 4.2.  Suppose that $f$ is critical.  Then with the ordered basis of $H^{1,1}(X)$: $H$, $E_1$, ${\cal S}_{0,3}$, $\dots$, ${\cal S}_{0,k}$, ${\cal P}_{0,1}$, $\dots$, ${\cal P}_{0,k-2}$ $E_k$, the action on cohomology classes $f_X^*$ on $H^{1,1}(X)$ is given by
$$\eqalign{f_X^*\ :\  & E_1 \mapsto {\cal S}_{0,3}\mapsto \cdots \mapsto {\cal S}_{0,k}\cr
& {\cal S}_{0,k}\mapsto \{\Sigma_0\} = H-E_1-{\cal S}_{0,3} - \cdots-{\cal S}_{0,k}-{\cal P}_{0,1}-\cdots-{\cal P}_{0,k-2}-E_k\cr
& {\cal P}_{0,1} \mapsto {\cal P}_{0,2}\mapsto \cdots\mapsto {\cal P}_{0,k-2} \mapsto E_k\cr
&E_k \mapsto \{\Sigma_\beta\}=H-{\cal P}_{0,1}-\cdots-{\cal P}_{0,k-2}-E_k\cr
& H \mapsto 2H - E_1-{\cal P}_{0,1}-\cdots-{\cal P}_{0,k-2}-E_k.\cr}\eqno{(4.4)}$$
The action on cohomology ${f_X^{-1}}^*$ is similar. In fact, the matrix representation for $f_X^*$ and ${f_X^{-1}}^*$ are the same up to the order of basis. Furthermore the spectral radius is given by the largest root of $x^{2k-1} - (x^k-1)/(x-1)$.

\medskip

Let us denote $ \Delta^c_k$ the largest root of the polynomial $x^{2k-1} - (x^k-1)/(x-1)$. Using the fact that the dynamical degree is upper semi-continuous, we obtain the largest dynamical degree in critical cases. 

\proclaim Theorem 4.3.   In the critical case, $\delta(f_X)\le \Delta_k^c$ and  $\delta(f_X^{-1})\le \Delta_k^c$. 

\noindent{\bf \S 5. Periodic Mappings }   
In this section, we consider the induced birational mapping $f_X$ such that the orbit of $\Sigma_\gamma$ ends up with $\Sigma_{\beta \gamma}$, that is  for some $n_\star \ge 0$
$$f_X : \Sigma_\gamma \mapsto \Sigma_{BC} \mapsto f_X\Sigma_{BC} \mapsto \cdots \mapsto f^{n_\star}\Sigma_{BC} = \Sigma_{\beta\gamma}$$ and $f^j_X \Sigma_{BC} \not\in \Sigma_0 \cup \Sigma_\beta \cup \Sigma_\gamma$ for $j=0, \dots, n_\star-1$.  Let $\pi_Z: Z\to X$ be a complex manifold obtained by blowing up the orbit of $\Sigma_{BC}=f_X\Sigma_\gamma$ and denote ${\cal F}_j$ the exceptional divisor over $f^{j-1}_X\Sigma_{BC}$ for $j=1, \dots, n_\star+1$. 

\proclaim Lemma 5.1. If there exists a positive integer $n_\star$ such that $f_X^{n_\star} \Sigma_{BC} = \Sigma_{\beta\gamma}$, the dynamical degree is given by the largest root of the polynomial
$$\chi_{k,n_{\star}}(x)= (x^{1+2 k+n_\star}-x^{2k+n_\star}-x^{1+k+n_\star}+x^{1+n_\star}+x^{2 k} -x^k-x+1)/(x-1).\eqno{(5.1)}$$

\noindent{\it Proof.} Since $f_X$ is well defined on $\Sigma_{BC}, \dots, f^{n_\star-1}\Sigma$, it suffices to check the mapping on $\Sigma_{\gamma}$ and $\Sigma_{\beta\gamma}$. By the induced map $f_Z$ the generic point on $\Sigma_\gamma$ map to a point on ${\cal F}_1$ :
$$\eqalign{f_Z: \Sigma_\gamma \ni  &[ x_0:x_1: \cdots: x_{k-1}: -\alpha_0 x_0 -\alpha_2 x_2 - \cdots - \alpha_{k-1} x_{k-1}] \cr&\mapsto (x_2/x_0, \dots, x_{k-1}/x_0, 0, x_0/(\beta_0 x_0+ x_1))_{f_1} \in {\cal F}_1}$$
where we use a local coordinates near ${\cal F}_1$ : 
$$\pi_{f_1}: (x_1, \dots, x_{k-2}, s, \xi_k)_{f_1} \mapsto [ 1: x_1: \cdots: x_{k-2}: - \alpha_0 - \alpha_2 x_1 - \cdots - \alpha_{k-1} x_{k-2} + s: s \xi_k] $$
Also under the inverse map $f^{-1}_Z$ we have 
$$\eqalign{f^{-1}_Z: &\Sigma_C \ni [ x_0 :\cdots: x_{k-2}: -\alpha_0 x_0 - \alpha_2 x_1- \cdots - \alpha_{k-1} x_{k-2} : x_k] \cr & \mapsto (x_0,0,x_2,\dots, x_{k-2}, x_k/x_0)_{f_{n_\star}} \in {\cal F}_{n_\star}}$$ 
where we use a local coordinates near ${\cal F}_{n_\star}$ : 
$$\pi_{f_{n_\star}}: (x_0,s,x_2, \dots, x_{k-2}, \xi_k)_{f_{n_\star}} \mapsto [ x_0: -\beta_0 x_0 + s: x_2:\cdots: x_{k-2}: - \alpha_0 - \alpha_2 x_2 - \cdots - \alpha_{k-1} x_{k-1} : s \xi_k]. $$
It follows that $f_Z$ is a local diffeomorphism at points on  $\Sigma_\gamma \cup \bigcup_{j=1}^{n_\star} {\cal F}_j$. Further more $f_Z$ doesn't have any exceptional hypersurfaces and therefore $f_Z$ is $1$-regular. To compute the action on $H^{1,1}(Z)$ let us choose the ordered basis $H$, $E_1$, ${\cal S}_{0,3}$, $\dots$, ${\cal S}_{0,k}$, ${\cal P}_{0,1}$, $\dots$, ${\cal P}_{0,k-2}$ $E_k$, ${\cal F}_{n_\star}$, $\dots$, ${\cal F}_1$. The action on cohomology classes $f_Z^*$ on $H^{1,1}(Z)$ is given by
$$\eqalign{f_X^*\ :\  & E_1 \mapsto {\cal S}_{0,3}\mapsto \cdots \mapsto {\cal S}_{0,k}\cr
& {\cal S}_{0,k}\mapsto \{\Sigma_0\} = H-E_1-{\cal S}_{0,3} - \cdots-{\cal S}_{0,k}-{\cal P}_{0,1}-\cdots-{\cal P}_{0,k-2}-E_k\cr
& {\cal P}_{0,1} \mapsto {\cal P}_{0,2}\mapsto \cdots\mapsto {\cal P}_{0,k-2} \mapsto E_k\cr
&E_k \mapsto \{\Sigma_\beta\}=H-{\cal P}_{0,1}-\cdots-{\cal P}_{0,k-2}-E_k\cr
&{\cal F}_{n_\star} \mapsto \cdots \mapsto {\cal F}_1  \mapsto \{ \Sigma_\gamma\} = H- E_1 - {\cal F}_{n_\star}\cr
& H \mapsto 2H - E_1-{\cal P}_{0,1}-\cdots-{\cal P}_{0,k-2}-E_k.\cr}\eqno{(5.2)}$$
The spectral radius of the action given by $(5.2)$ is the largest root of $\chi_{k,n_{\star}}(x)$ \qed

\proclaim Lemma 5.2.  If $n_\star > (k^2+k)/(k-1)$ then $f$ has exponential degree growth.

\noindent{\it Proof.} The derivative of $\chi$ at $x=1$ is negative if $n_\star > (k^2+k)/(k-1)$. It follows that $n_\star > (k^2+k)/(k-1)$ implies that $\chi$ has a real root which is strictly bigger than $1$.\qed

\proclaim Lemma 5.3.  For a critical map, $n_\star\ge k-1$.

\noindent{\it Proof.} Since we have $\Sigma_{BC} \subset \Sigma_k$, $\Sigma_{\beta\gamma} \not\subset \Sigma_j$, for $j\ge 2$, and $f: \Sigma_k \mapsto \Sigma_{k-1} \mapsto \cdots\mapsto \Sigma_1$, it requires at least $k-1$ iterations for $\Sigma_{BC}$ to be mapped to $\Sigma_{\beta\gamma}$.  \qed

\proclaim Lemma 5.4. If $k>3$ and $n_\star>k+2$ then the dynamical degree is strictly bigger than $1$.  If $k=3$, the dynamical degree for $n_\star=6=(k^2+k)/(k-1)$  is equal to $1$ and the dynamical degree for $n_\star \ge 7$ is strictly bigger than $1$.

\noindent{\it Proof.} The second derivative of $(x-1) \chi_{k, n_\star}(x)$ at $x=1$ is $2 ((1-k) n_\star + k (k+1))$. It follows that $  \chi_{k, n_\star}'(1) <0$ and therefore the dynamical degree is strictly bigger than $1$  if and only if $n_\star > (k^2+k)/(k-1)$. Since $k+3 = (k^2 + k + (k-3))/(k-1)$, $k+3 > (k^2+k)/(k-1)$ if $k >3$ and $k+3 = (k^2+k)/(k-1)$ if $k =3$. \qed

In case $n_\star=k+2$, we have $\chi'_{k, n_\star}(1)=2$. Using the computer we checked for each $4\le k\le 20$ and found that $\chi_{k, k+2}(x)$ has complex conjugate pair or roots whose modulus is strictly bigger than $1$.

%


\proclaim Theorem 5.5. Let $f$ be a critical map with $k\ge3$.  If there exists a positive integer $n_\star$ such that $f^{n_\star} \Sigma_{BC} = \Sigma_{\beta\gamma}$, then one of the following must occur:
$$ \eqalign{ (i) & {\rm \ \ If\ }n_\star =k-1 {\rm \ then\ the \ mapping\ is\ periodic\ with\ period\ }3k-1\cr
(ii) & {\rm \ \ If\ }n_\star =k\phantom{-1} {\rm \ then\ the \ mapping\ is\ periodic\ with\ period\ }4k\cr
(iii) & {\rm \ \ If\ }n_\star =k+1 {\rm \ then\ the \ mapping\ is\ periodic\ with\ period\ }3k(k+1).\cr}$$

\noindent{\it Proof.} Using $(5.1)$ it is not hard to show that $\chi_{k, k-1} = (x^{3k-1}-1)$, $\chi_{k,k} = (x^k-1) (x^{2k}+1)$ and $\chi_{k, k+1} = (x^{k+1}-1)(x^{2k}-x^k+1)$. In each case $f$ is linear fractional. Suppose $f[x_0: x_1:\cdots:x_k] = [\sum_i a_{0,i} x_i:\sum_i a_{1,i} x_i:\cdots:\sum_i a_{k,i} x_i]$. Using the fact that $E_1$ and $E_k$ is fixed under $f_X$ we see that $a_{j,1} =0$ for all $j\ne 1$ and $a_{j,k} =0$ for all $j \ne k$. Since we are in the projective space we may assume that $a_{1,1}=a_{k,k}=1$. For each fixed co-dimension $j$ subspace we obtain $j+1$ equations on $a_{j,i}$. We continue this procedure for other fixed linear subspaces to conclude that the mapping is actually periodic.  \qed

The case of dimension $k=2$ is not covered by Theorem 5.5; the numbers corresponding to the cases {\it (i--iii)} are 5, 8, and 18.  These are all found to occur in [BK1], where it was shown that there are also the possibilities of period 6, 12, and 30.   If $k=2$ then  $\chi_{2, n_\star}'(1)=6(6-n_\star)$ and therefore we have more possibilities for periodic mappings, that is $n_\star$ could be $k+2$ and $k+3$ which correspond to the cases of period 12 and period 30. The mapping with period 6 occurs when $\Sigma_{BC} = \Sigma_{\beta\gamma}$.  This cannot happen in dimension $3$ or higher, since $\Sigma_{BC}$ and $\Sigma_{\beta\gamma}$ are linear spaces of positive dimension, and there exists a point $[1:1-\beta_0: x_2: \cdots: x_{k-2}: -\alpha_0 - \alpha_3 x_2 - \cdots -\alpha_{k-1} x_{k-2}:0] \in \Sigma_{BC} \setminus \Sigma_{\beta\gamma}$. In the case of dimension $k=3$,  [BK2] shows that the only possible periods are 8 and 12 (which correspond to cases {\it (i)} and {\it (ii)} in Theorem 5.5); the possibility  $n_\star=k+1=4$ does not occur in dimension 3. 

\proclaim Theorem 5.6. If $\alpha_{k-1}=(-1)^{1/k}$ and 
$$\beta=(\alpha_{k-1}^{k-1},1,0,\dots,0)\ \ {\rm and\ \ }\alpha=(\alpha_{k-1}^{k-2}/(1-\alpha_{k-1}),0,\alpha_{k-1}^{k-2},\dots, \alpha_{k-1}^2,\alpha_{k-1},1)\eqno{(5.3)}$$
then $f_{\alpha,\beta}$ is periodic with period $4k$.

\noindent{\it Proof.} It is suffices to show that with these choices of parameter values we have $f^{k} \Sigma_{BC} = \Sigma_{\beta\gamma}$. Let us set ${\cal A}:=-(\alpha_0 x_0 + \alpha_2 x_1+ \cdots+ \alpha_{k-1} x_{k-2})$. The generic point $p \in \Sigma_{BC}$ can be written as $[x_0:x_1:\cdots:x_{k-2}:{\cal A}:0]$. The last coordinate of $f(p)$ is given by $$x_0 ( \alpha_0 x_0 + \alpha_2 x_2+ \cdots + \alpha_{k-2} x_{k-2} - \alpha_{k-1} {\cal A}).$$ Since $\alpha_{k-1} \alpha_j = \alpha_{j-1}$ for $j=3, \dots, k$ and $\alpha_0 - \alpha_{k-1} \alpha_0 = \alpha_{k-1}^{k-2}$, the last coordinate of $f(p)$ becomes $$x_0( \alpha_{k-1}^{k-2} x_0 - \alpha_{k-1}^{k-1}x_1) = -\alpha_{k-1}^{k-1} ( - \alpha_{k-1}^{-1} x_0 + x_1) x_0= -\alpha_{k-1}^{k-1}  x_0\beta \cdot x. $$ It follows that $f: \Sigma_{BC} \ni [x_0:x_1:\cdots:x_{k-2}:{\cal A}:0] \mapsto [ x_0:x_2: \cdots: x_{k-2}: {\cal A}:0:-\alpha_{k-1}^{k-1} x_0]$ and using that $\alpha_{k-1}^{k} = -1$ we have 
$$f: \Sigma_{BC} \mapsto \{x_{k-1} = 0 \} \cap \{x_0 - \alpha_{k-1} x_k = 0 \}.$$
From $(2.1)$ it is not hard to see that $$\eqalign{& f:  \{x_{k-1} = 0 \} \mapsto  \{x_{k-2} = 0 \} \mapsto \cdots \mapsto  \{x_{1} = 0 \} \cr & f: \{x_0 - \alpha_{k-1} x_k = 0 \}\mapsto \{x_0 - \alpha_{k-1} x_{k-1} = 0 \} \mapsto \cdots \mapsto \{x_0 - \alpha_{k-1} x_{2} = 0 \}. \cr}$$
That is $f^{k-1}\Sigma_{BC} =  \{x_{1} = 0 \} \cap \{x_0 - \alpha_{k-1} x_{2} = 0 \}$. Now let us map forward a point $p=[\alpha_{k-1} x_{2}:0:x_2:\cdots:x_k]$. Since $\beta\cdot p = - x_2$ we have
$$\eqalign{f: &[\alpha_{k-1} x_{2}:0:x_2:\cdots:x_k] \mapsto \cr &[ -\alpha_{k-1} x_2: -x_2:\cdots:-x_k: \alpha_{k-1} (\alpha_0 \alpha_{k-1} + \alpha_2) x_2+ \alpha_{k-1} \alpha_3 x_3 + \dots + \alpha_{k-1} x_k].}$$
It follows that $\beta\cdot fp = - \alpha_{k-1}^k x_2 -x_2 =0$ and $\gamma\cdot fp = \alpha_{k-1} (\alpha_0(\alpha_{k-1}-1)+ \alpha_2) x_2 + (-\alpha_2+ \alpha_{k-1} \alpha_3) x_3 + \cdots + (-\alpha_{k-1}+ \alpha_{k-1}) x_k =0$ and therefore $f^k \Sigma_{BC} = \Sigma_{\beta\gamma}$.\qed

\bigskip

\noindent{\bf \S 6. Non-periodic maps; integrability. } 
Let us consider the critical map given by $\alpha= (a, 0,1,\dots, 1)$ and $\beta= (0,1,0,\dots,0)$:
$$\eqalign{f[x_0:\cdots:x_k] &= [x_0 x_1: x_2 x_1: \cdots : x_k x_1: x_0 (a x_0 + x_2 + \cdots +x_k)]\cr f^{-1}[x_0:\cdots:x_k]&=[x_0 x_k: x_0(a x_0+x_1+\cdots+x_{k-1}):x_1 x_k: \cdots: x_{k-1} x_k].}\eqno{(6.1)}$$
It follows that we have
$$f: \ \ \ \eqalign{&\Sigma_\beta \mapsto e_k \rightsquigarrow \Sigma_{0\,1\,\dots\,k-2}\rightsquigarrow \Sigma_{0\,1\,\dots
\,k-3}\rightsquigarrow\cdots \rightsquigarrow \Sigma_{0\,1} = \Sigma_{0\beta} \rightsquigarrow \Sigma_0\cr
& \Sigma_0 \mapsto \Sigma_{0k}\mapsto \Sigma_{0\,k-1\,k}\mapsto \cdots \mapsto \Sigma_{0\,3\,\dots\,k} \mapsto e_1 \rightsquigarrow \Sigma_B \cr
}\eqno{(6.2)}$$
In addition, since $\beta_0=0$ we  have 
$$f:\Sigma_B=\Sigma_k \mapsto\Sigma_{k-1} \mapsto \cdots \mapsto \Sigma_1=\Sigma_\beta  \eqno(6.3)$$ 
so we expect to find (after blowing up) a closed orbit of hypersurfaces containing  $\Sigma_\beta$, $\Sigma_0$, and $\Sigma_B$.

Our first task will be to make $f$ 1-regular.
For $j=1, \dots, k-1$ let us set $q_j = [0: \cdots:0: 1:-1:0: \cdots :0]$, the point whose $j$-th coordinate is $1$, whose $j+1$-th coordinate is $-1$, and every other coordinate is zero. Let us consider a complex manifold $\pi_1: Z_1 \to \P^k$ obtained by blowing up the $k+1$ points $e_1$, $e_k$, and  $q_j$, $j=1, \dots, k-1$.  We denote by ${\cal Q}_j$ the exceptional divisor over the point $q_j$.  We also denote by $E_1$ and $E_k$ the exceptional divisors over the points $e_1$ and $e_k$.
\proclaim Lemma 6.1. The induced map $f_{Z_1}$ is a local diffeomorphism at generic points of ${\cal Q}_j, j=2, \dots, k-1$. Furthermore we have dominant maps
$$f_{Z_1} : {\cal Q}_{k-1}\mapsto {\cal Q}_{k-2} \mapsto \cdots \mapsto {\cal Q}_2 \mapsto {\cal Q}_1$$

\noindent{\it Proof.} Let us consider the local coordinates near ${\cal Q}_{k-1}$ and ${\cal Q}_{k-2}$ 
$$\eqalign{&\pi_{k-1}: Z_1\ni(s, \xi_1, \cdots, \xi_{k-1})_{k-1} \mapsto [s: s \xi_1:\cdots: s \xi_{k-2}: 1+s \xi_{k-1}:-1] \in \P^k\cr
&\pi_{k-2}: Z_1\ni(s, \xi_1, \cdots, \xi_{k-2},\xi_k)_{k-2} \mapsto [s: s \xi_1:\cdots: s \xi_{k-3}: 1+s \xi_{k-2}:-1: s \xi_k] \in \P^k.}$$
Note that in those coordinates, $\{(s, \xi_1, \cdots, \xi_{k-1})_{k-1}: s=0\}= {\cal Q}_{k-1}$, and  we see that 
$$ f_{Z_1} :{\cal Q}_{k-1} \ni(0, \xi_1, \cdots, \xi_{k-1})_{k-1} \mapsto (0,\xi_2,\dots, \xi_{k-1}, {a \over \xi_1}+\xi_2+\cdots \xi_{k-1})_{k-2}\in{\cal Q}_{k-2}.$$
It follows that $f_{Z_1}$ is locally diffeomorphic at generic points of ${\cal Q}_{k-1}$. For $j=2,\dots,k-2$, the proof is identical.  \qed

By constructing $Z_1$, we create three new exceptional hypersufaces including $E_1$ and $E_k$ for each $f_{Z_1}$ and $f_{Z_1}^{-1}$. Let us consider the local coordinates $\pi_1: (s, \xi_1,\xi_3 \dots, \xi_k)_1 \mapsto [s: 1+s \xi_1: -1: s \xi_3: \cdots: s\xi_k]$ near ${\cal Q}_1$. We also use the local coordinates $\pi_{e_1} :(s, \xi_2, \dots, \xi_k)_{e_1} \mapsto [s: 1: s \xi_2: \cdots:s \xi_k]$ near $E_1$. With these coordinates, we see that
$$f_{Z_1}: {\cal Q}_1 \ni (0, \xi_1,\xi_3 \dots, \xi_k)_1\mapsto (0,\xi_3, \dots, \xi_k,-1)_{e_1} \in E_1 \cap \{x_0+x_k=0\}\subset E_1.$$
Similarly with the local coordinates $\pi_{e_k} :(s, \xi_1, \dots, \xi_{k-1})_{e_k} \mapsto [s:  s \xi_1: \cdots:s \xi_{k-1}:1]$ near $E_k$ and the local coordinates near ${\cal Q}_{k-1}$ defined above, we have
$$f^{-1}_{Z_1}: {\cal Q}_{k-1} \ni (0, \xi_1,\xi_2 \dots, \xi_{k-1})_1\mapsto (0,-1,\xi_1, \dots, \xi_{k-2})_{e_k} =E_k \cap \{x_0+x_1=0\}\in E_k.$$
Thus we have
$$f_{Z_1} : E_k \cap \{x_0+x_1=0\} \rightsquigarrow  {\cal Q}_{k-1} \qquad {\rm and} \qquad f^{-1}_{Z_1}:E_1 \cap \{x_0+x_k=0\}\rightsquigarrow {\cal Q}_1\eqno{(6.4)}$$

\proclaim Lemma 6.2. $f_{Z_1}^{k} \Sigma_{BC}=E_k \cap \{x_0+x_1=0\}$ and $f_{Z_1}^{-k} \Sigma_{\beta\gamma} =E_1 \cap \{x_0+x_k=0\}$.

\noindent{\it Proof.} Let us consider the forward map first. A generic point $p$ in $\Sigma_{BC}$ can be written as $[x_0: \cdots: x_{k-2}: -a x_0-x_1-\cdots-x_{k-2}:0]$. Using $(6.1)$ we have 
$$f_X: [x_0: \cdots: x_{k-2}: -a x_0-x_1-\cdots-x_{k-2}:0] \mapsto [ x_0 x_1: x_2 x_1: \cdots: (-a x_0-x_1-\cdots-x_{k-2})x_1: 0:-x_0x_1].$$
It follows that $f_X\Sigma_{BC} \subset \{ x_{k-1}=0, x_0+x_k=0\}$. Since $\Sigma_{BC}$ is not indeterminate for $f$ we see that $f_X\Sigma_{BC} = \{ x_{k-1}=0, x_0+x_k=0\}$. Note that $f^{k-2} \{x_{k-1}=0\} = \{x_1=0\}$,  $f^{k-2} \{x_0 + x_k=0\} = \{ x_0+x_2=0\}$ and therefore $f_X^{k-1}\Sigma_{BC} = \{ x_{1}=0, x_0+x_2=0\}$
Using the coordinates near $E_k$ we see that 
$$f_X :  \Sigma_1 \ni [x_0 :0: x_2 : \cdots: x_k] \mapsto (0, {x_2 \over x_0}, \dots,{x_k\over x_0})_{e_k} \in E_k$$
Thus we have 
$$f_X : f_X^{k-1}\Sigma_{BC}  \mapsto (0,-1,\xi_2, \dots, \xi_{k-1})_{e_k} \in E_k .$$
The argument for $f^{-1}$ is essentially identical. \qed

Now let us construct a complex manifold $\pi_2: Z_2 \to Z_1$ obtained by blowing up the sets $f_X^j \Sigma_{BC}, j=0, \dots, k$, ${\cal Q}_j, j=1,\dots k-1$ and $f^{-j}_X \Sigma_{\beta\gamma}, j=k,\dots, 0$. We denote ${\cal F}_j$ the exceptional divisor over the set $f_X^j \Sigma_{BC}$, and we also denote ${\cal H}_j$ the exceptional divisor  over the set $f^{-j}_X \Sigma_{\beta\gamma}$.

\proclaim Lemma 6.3. The induced map $f_{Z_2}$ is a local diffeomorphism at a generic points of $\bigcup_j {\cal F}_j \cup \bigcup_j{\cal H}_j$. Thus $f_{Z_2}$ has four exceptional hypersurfaces $\Sigma_0$, $\Sigma_\beta$, $E_1$, and $E_k$. 

\noindent{\it Proof.} It suffices to check at the points in $\Sigma_\gamma \cup {\cal F}_k \cup {\cal Q}_1 \cup {\cal H}_0$. Let us define local coordinates :
$$\eqalign{& (x_1,\dots, x_{k-2},s, \xi) \mapsto [1:x_1:\cdots:x_{k-2}: s-a-x_1-\cdots-x_{k-2}: s \xi]\qquad {\rm near\ } {\cal F}_1\cr
& (s, \eta, \xi_2, \dots, \xi_{k-1}) \mapsto (s, -1+ s \eta, \xi_2, \dots, \xi_{k-1})_{e_k}\qquad{\rm near \ } {\cal F}_{k}\cr
& (s, \xi_2, \dots, \xi_{k-2},\eta) \mapsto (s, \xi_2, \dots, \xi_{k-2}, -1+ s \eta,)_{e_1}\qquad{\rm near \ } {\cal H}_{k}\cr 
& (s, x_2,\dots, x_{k-1}, \xi) \mapsto [1:s:x_2:\cdots:x_{k-1}: -a-x_2-\cdots-x_{k-1}+ s \xi]\qquad {\rm near\ } {\cal H}_1.\cr}$$ 
Then we have
$$f_{Z_2}:\ \ \eqalign{& \Sigma_\gamma \ni [x_0: x_1: \cdots: x_{k-1}: -a x_0 -x_2-\cdots-x_{k-1}] \mapsto (x_2, \dots, x_{k-1}, 0, {x_0 x_1^{-1}}) \in {\cal F}_1 \cr
&{\cal F}_k \ni  (0,\eta, \xi_2, \dots, \xi_{k-1}) \mapsto (0, \xi_2, \dots, \xi_{k-1}, -\eta+a +\xi_2+ \cdots+\xi_{k-1}) \in {\cal Q}_{k-1}\cr 
&{\cal Q}_1\ni (0,\xi_1, \xi_3, \dots, \xi_k) \mapsto (0, \xi_3, \dots, \xi_k, a+\xi_1+\xi_3+\cdots+ \xi_k) \in {\cal H}_k\cr
&{\cal H}_1 \ni (0,x_2,\dots, x_{k-1}, \xi) \mapsto [1: x_2: \cdots:x_{k-1}: -a-x_2-\cdots-x_{k-1}:\xi] \in \Sigma_C.}$$
It follows that the induced map $f_{Z_2}$ is local diffeomorphism on the orbit of $\Sigma_\gamma$. \qed

For other two exceptional hypersurfaces $\Sigma_0$ and $\Sigma_\beta$, we construct a blowup space $\pi_{\cal Z} : {\cal Z} \to Z_2$ obtained by blowing up the strict transform of the sets $\Sigma_{0\,k},\dots, \Sigma_{0\,3\,\dots\,k-1\,k}, \Sigma_{0\,\dots\, k-2},\dots, \Sigma_{01}$ in $(6.4)$. We use the same notation for the space $X$ in \S3. That is, the exceptional divisors over  $\Sigma_{0\,1\,\dots\, k-2},\dots,\Sigma_{0\,1}$ are ${\cal P}_{0,k-2},\dots,{\cal P}_{0,1}$ and ${\cal S}_{0,j}$ are the exceptional divisors over $\Sigma_{0\,j\,\dots\,k}$ for all $j \ge 3$. Since are only consider a generic point on these exceptional divisors, the same computations as in \S3 and \S4 work and thus we conclude that $f_{\cal Z}$  is local diffeomorphic at a generic points on these new exceptional divisors as well as $E_1,E_k,\Sigma_0$ and $\Sigma_\beta$. It follows that: 
\proclaim Lemma 6.4. The induced map $f_{\cal Z}$ has no exceptional hypersurface and therefore $f_{\cal Z} $ is $1$-regular.

Now to compute the dynamical degree we use the following basis of $Pic({\cal Z})$:
$$H, E_1, {\cal S}_{0,3}, \dots, {\cal S}_{0,k}, {\cal P}_{0,1}, \dots, {\cal P}_{0,k-2}, E_k, {\cal H}_{0}, \dots, {\cal H}_k, {\cal Q}_{k-1},\dots, {\cal Q}_{1},{\cal F}_k,\dots, {\cal F}_0. $$ 
Using $(6.2)$, Lemma 6.1 and Lemma 6.2 we have:

\proclaim Lemma 6.5. The action on cohomology $f_{\cal Z}^*$ is given by 
$$\eqalign{f_{\cal Z}^*\ :\  & E_1 \mapsto {\cal S}_{0,3}\mapsto \cdots \mapsto {\cal S}_{0,k}\mapsto \{\Sigma_0\} \cr
& {\cal P}_{0,1} \mapsto {\cal P}_{0,2}\mapsto \cdots\mapsto {\cal P}_{0,k-2} \mapsto E_k \mapsto \{\Sigma_\beta\}\cr
&{\cal H}_{3k+1} \mapsto {\cal H}_{3k} \mapsto \cdots \mapsto {\cal H}_1 \mapsto \{\Sigma_\gamma\}\cr
& H \mapsto 2H - E_1-{\cal P}_{0,1}-\cdots-{\cal P}_{0,k-2}-E_k-{\cal H}_0-{\cal H}_{k}-\sum_{j=2}^{k-1} {\cal Q}_j-{\cal F}_k.\cr}\eqno{(6.5)}$$
where
$$\eqalign{ &  \{\Sigma_0\} = H-E_1-{\cal S}_{0,3} - \cdots-{\cal S}_{0,k}-{\cal P}_{0,1}-\cdots-{\cal P}_{0,k-2}-E_k-{\cal H}_k-\sum_{j=1}^{k-1} {\cal Q}_j - {\cal F}_k\cr &  \{\Sigma_\beta\}=H-{\cal P}_{0,1}-\cdots-{\cal P}_{0,k-2}-E_k-{\cal H}_0 - \sum_{j=2}^{k-1} {\cal Q}_j-{\cal F}_k-{\cal F}_{k-1} \cr &\{\Sigma_\gamma\} = H-E_1-{\cal H}_0-{\cal H}_{k}-\sum_{j=2}^{k-1} {\cal Q}_j}$$

\proclaim Theorem 6.6. For every $k>3$, the map $f$ defined in $(6.1)$ has quadratic degree growth. 

\noindent{\it Proof.} The characteristic polynomial of the action on cohomology given in $(6.5)$ is given by
$$ \hat\chi_k (x) = \pm(x^k-1)(x^{k+1}-1)(x^{3k-1}-1).$$
It follows that $1$ is a zero of $ \hat\chi_k (x)$ with multiplicity $3$. Furthermore there is unique (up to scalar multiple) eigenvector $v$ corresponding to an eigenvalue $1$:
$$\eqalign{v\ =\  &-(k+1) H + (k-1) E_1 + \sum_{j=1}^{k-2} j \,{\cal S}_{0,k+1-j} + \sum_{j=1}^{k-2} j \,{\cal P}_{0,j} + (k-1) E_k\cr & + \sum_{j=0}^{k-1} {\cal H}_j + k {\cal H}_k+ (k-1) \sum_{j-1}^{k-1} {\cal Q}_j + k {\cal F}_k + \sum_{j=0}^{k-1}{\cal F}_k}$$
It follows that the Jordan decomposition has $3\times 3$ block with 1 on the diagonal.  Thus the powers of this matrix grow quadratically.  \qed

We say that a rational function $\varphi$ is an {\it integral} of $f$ if $\varphi=\varphi\circ f$ at generic points.  Some integrals of $f$ have been found  (see [KLR], [KL] [CGM1], [GKI]).  Here we describe a rationale for re-finding these (known) integrals.  We start by finding homogeneous polynomials $p$ which are invariant in the sense that 
$$p\circ f = J \cdot p   \eqno(6.6)$$
 where  $J =x_0 (\beta \cdot x)^{k-1} (\gamma \cdot x) $ is the Jacobian of  $f$.  This is the same as finding a meromorphic $k$-form $\eta$, written as ${dx_1\wedge\cdots\wedge dx_k / p(1,x_1,\dots,x_k) }$ on the affine coordinate chart $\{x_0=1\}$, and which is invariant in the sense that $f^*\eta=\eta$.  If $p_1,\dots, p_r$ satisfy (6.6), then $\sum\lambda_jp_j$ will also satisfy (6.6).  And the quotient of any two of these polynomials  $\sum\lambda_jp_j$ will give an integral.

Since $f$ has degree 2,  $J$ has degree $k+1$, so we look for polynomials $p$ of degree $k+1$, so that the degrees of $p\circ f$ and $J p$ will both be $2(k+1)$.  The invariant rational functions will then be given as quotients of invariant functions $h=p_1/p_2$.  Recall that $f$ maps $\Sigma_\beta$ to $e_k$.  Thus by (6.6) we see that $p$ will vanish to order at least $k-1$ at $e_k$, since $J$ vanishes to order $k-1$ at $\Sigma_\beta$.  Similarly, since $f(\Sigma_0)=\Sigma_{0,k}$, we see that $p$ must vanish at $\Sigma_{0,k}$.  Now starting with a point $z\in\Sigma_{0,k}$, we have $f(z)\in\Sigma_{0,k-1,k}$, so by (6.6), $p$ vanishes to order at least 2 on $\Sigma_{0,k-1,k}$.  Continuing this way, we see that $p$ vanishes to order at least $k-j$ on $\Sigma_{j+1,j+2,\dots,k}$ for $1\le j\le k-1$.  Finally, since $J$ vanishes on $\Sigma_\gamma$, and $f(\Sigma_\gamma)=\Sigma_{BC}$, we see that $p$ vanishes on $\Sigma_{BC}$.  Iterating this, we see that $p$ must vanish on $f^j(\Sigma_{BC})$ for $j=0,\dots,k$.

In order to refine the equation (6.6), we define $A=\alpha\cdot x$, $B=\beta\cdot x = x_1$ and $C=x_0$, so $J=AB^{k-1}C$.  Looking for linear functions which vanish on certain of the sets above, we define three families: 
$$\ell_j(x):=x_j, \ \ 0\le j\le k,\ \ \   m_j(x):=x_0+x_j, \ \ 1\le j\le k, \ \ \ n_j:=x_0 + x_j+x_{j+1}, \ \ 1\le j\le k-1,$$
and $n_0=n_{k-1}\circ f$ and $m_0 = A + x_1 = a x_0 + x_1 + \cdots + x_k$.   Thus we have a refined form of (6.6)
$$\eqalign{ &\ell_j\circ f=B\ell_{j+1}, \ \ \ m_j\circ f = B m_{j+1}, \ \ 1\le j\le k-1, \ \ \ n_j\circ f=B n_{j+1}, 1\le j\le k-2\cr
&\ell_k\circ f = A \ell_0, \  \ \ \ell_0\circ f = C\ell_1, \ \ \ m_k\circ f = C m_0,  \ \ \ m_0\circ f = A m_1,  \ \ \ n_0\circ f = ABC n_1.}\eqno(6.7)$$

By (6.7) it is evident that $p_0:=\ell_0\ell_1\cdots\ell_k$ and $p_1:=m_0 m_1\cdots m_k$ satisfy (6.6).  If $k\ge 3$, then  $p_2:= n_0\cdots n_{k-1}$ also satisfies (6.6).  

Now let us use the notation ${\bf j}$ for the product $\ell_jm_j$.  If $k\ge 5$ is odd, we define 
$$\Phi_{\rm even}:={\bf 0\ 2 \ 4\cdots (k-1)}\ \ \ \ \ \Phi_{\rm odd}:={\bf 1\ 3\ 5\cdots k}.$$
If $q$ is a polynomial for which $q\circ f$ is divisible by $J$, we let $T$ denote the operator $T(q) = q\circ f\cdot J^{-1}$, so that (6.6) holds exactly when $p$ is a fixed point of $T$.   
By (6.7) we have ${\bf j}\circ f={\bf (j+1)} \ B^2$ for $1\le j\le k-1$; and ${\bf k}\circ f={\bf 0}\  AC$, and ${\bf 0}\circ f={\bf 1}\  AC$.
Thus  $T\Phi_{\rm even}=  \Phi_{\rm odd}$, and $T\Phi_{\rm odd} = \Phi_{\rm even}$.   We conclude that 
$p_3:=\Phi_{\rm even}+\Phi_{\rm odd}$ satisfies (6.6).

If $k>5$ is even, we consider two functions:
$$\Psi_a: = {\bf 0}\  n_1\  {\bf  3\ 5\ 7\ \cdots (k-1)}, \ \ \ \ \  \Psi_b :={\bf 1\ 3\ 5\ 7\ \cdots (k-1)}\ \ell_k .$$
By  (6.7) we see that $\Psi_b\circ f = A\ \ell_0 \ {\bf 2\ 4\ 6\ 8\ \cdots\ k}\ B^k= J\ B\ {\bf 2\ 4\ 6\ 8\ \cdots\ k}$.  Thus $T\Psi_b=\ell_1\ {\bf 2\ 4\ 6\ 8\ \cdots\ k}$.  Applying $T$ to $\Psi_a$, we have 
$$\eqalign{T\Psi_a = {\bf 1}\ n_2\ &  {\bf 4\ 6\ 8\ \cdots\ k}, \ \ \ T^2\Psi_a = {\bf 0\ 2}\ n_3\ {\bf 5\ 7\ 9 \ \cdots\ (k-1)k}, \cr
\dots\ \ \ &T^{k-2}\Psi_a = {\bf 0\ 2\ 4\ \cdots (k-2)}\ n_{k-1}  }\eqno(6.8)$$
Now we claim that 
$$p_3 = \left(\Psi_a+T\Psi_a+\cdots  + T^{k-2}\Psi_a\right) + \left( \Psi_b + T\Psi_b\right)$$
satisfies (6.6).  For this, it suffices to have
$$\Psi_a+\Psi_b = T^{k-1}\Psi_a + T^2\Psi_b.\eqno(6.9)$$ 
Applying $T$ to $T^{k-2}\Psi_a$ in (6.8), and using (6.7), we find $T^{k-1}\Psi_a =  n_0\ {\bf 1\ 3\ 5\ \cdots\ (k-1)}$.  Now apply $T$ to the expression for  $T\Psi_b$ found above, we find $T^2\Psi_b={\bf 0}\ \ell_2\ {\bf 3\ 5\ \cdots\ (k-1)}$.   Thus (6.9) is a consequence of the simple identity
$${\bf 0}\ n_1 + {\bf 1}\ \ell_k    =  n_0\ {\bf 1} + {\bf 0}\ \ell_2,$$
and we conclude that $p_3$ satisfies (6.6).

\centerline{\bf References}

\medskip
\item{[BK1]}  E. Bedford and  Kyounghee Kim, Periodicities in linear fractional recurrences: degree growth of birational surface maps. 
Michigan Math. J. 54 (2006), no. 3, 647--670. 

\item{[BK2]}   E. Bedford and  Kyounghee Kim, Periodicities in 3-step linear fractional recurrences, preprint.

\item{[CL]}  E. Camouzis and  G. Ladas, {\sl  Dynamics of third-order rational difference equations with open problems and conjectures}. Advances in Discrete Mathematics and Applications, 5. Chapman \&\ Hall/CRC, Boca Raton, FL, 2008.

\item{[CGM1]} A.  Cima, A.  Gasull,  V. Ma\~nosa, Dynamics of some rational discrete dynamical systems via invariants. Internat. J. Bifur. Chaos Appl. Sci. Engrg. 16 (2006), no. 3, 631--645.

\item{[CGM2]} A.  Cima, A.  Gasull,  V. Ma\~nosa, Dynamics of the third order Lyness' difference equation. J. Difference Equ. Appl. 13 (2007), no. 10, 855--884.

\item{[CGM3]} A.  Cima, A.  Gasull,  V. Ma\~nosa, Some properties of the $k$-dimensional Lyness' map. J. Phys.\ A 41 (2008), no. 28, 285205, 18 pp.


\item{[CsLa]} M. Cs\"ornyei and M. Laczkovich, Some periodic and non-periodic recursions, Monatsh.\ Math.\ 132 (2001), 215--236.

\item{[dFE]}  T. de Fernex and L. Ein,  Resolution of indeterminacy of pairs.
{\sl Algebraic geometry}, 165--177, de Gruyter, Berlin, 2002.

\item{[DF]}  J. Diller and C. Favre, Dynamics of bimeromorphic maps of
surfaces, Amer. J. of Math., 123 (2001), 1135--1169.

\item{[GKI]} M. Gao,Y. Kato and M. Ito, Some invariants for $k$th-order Lyness equation. Appl. Math. Lett. 17 (2004), no. 10, 1183--1189.

\item{[GBM]}  L. Gardini,  G. Bischi, and C. Mira,  Invariant curves and focal points in a Lyness iterative process. Dynamical systems and functional equations (Murcia, 2000). Internat. J. Bifur. Chaos Appl. Sci. Engrg. 13 (2003), no. 7, 1841--1852.

\item{[G]}  M. Gizatullin, Rational $G$-surfaces. (Russian) Izv. Akad. Nauk SSSR Ser. Mat. 44 (1980), no. 1, 110--144, 239. 

\item{[GL]}  E. A. Grove and G. Ladas, {\sl Periodicities in nonlinear difference equations}, Adv. Discrete Math. Appl., 4, Chapman \&\ Hall and CRC Press, Boca Raton, FL, 2005.

\item{[H]} B. Hassett,  {\sl Introduction to Algebraic Geometry}, Cambridge U. Press, 2007.

\item{[HKY]}  R. Hirota, K. Kimura and  H. Yahagi,  How to find the conserved quantities of nonlinear discrete equations, J. Phys.\ A: Math.\ Gen.\ 34 (2001) 10377--10386.

\item{[KL]} V.I. Kocic and  G. Ladas,  {\sl Global Behaviour of Nonlinear Difference Equations of Higher Order with Applications}, Kluwer Academic Publishers 1993.

\item{[KLR]} V.I. Kocic, G. Ladas, and I.W. Rodrigues, On rational recursive sequences, J. Math. Anal. Appl 173 (1993), 127-157.

\item{[L]}  R. C. Lyness, Notes 1581, 1847, and 2952, Math. Gazette 26 (1942), 62; 29 (1945), 231; 45 (1961), 201.

\item{[M]}  C. T. McMullen,  Dynamics on blowups of the projective plane. Publ. Math. Inst. Hautes \'Etudes Sci. No. 105 (2007), 49--89.

\item{[Z]} C. Zeeman,  Geometric unfolding of a difference equation.  Lecture.

\bigskip

\rightline{E. Bedford:  bedford@indiana.edu}

\rightline {Department of Mathematics}

\rightline{Indiana University}

\rightline{Bloomington, IN 47405 USA}

\medskip\rightline{K. Kim: kim@math.fsu.edu}

\rightline {Department of Mathematics}

\rightline{Florida State University}

\rightline{Tallahassee, FL 32306 USA}

\bye